\documentclass[reqno]{amsart}

\usepackage{amssymb}
\usepackage{graphicx}
\usepackage{amscd}
\usepackage[pagebackref]{hyperref}
\usepackage{color}
\usepackage{tabularx}
\usepackage[table]{xcolor}
\usepackage{float}
\usepackage{graphics,amsmath,amssymb}
\usepackage{amsthm}
\usepackage{amsfonts}
\usepackage{latexsym}
\usepackage{epsf}
\usepackage{xifthen}
\usepackage{mathrsfs}
\usepackage{dsfont}
\usepackage{makecell}
\usepackage{subfig}
\usepackage{amsmath}
\usepackage{listings}
\usepackage{etoolbox}
\usepackage{fancyhdr}
\usepackage{pdflscape}
\usepackage[title,toc,titletoc]{appendix}
\usepackage{enumitem}
\usepackage[noadjust]{cite}
\usepackage{tikz}
\usetikzlibrary{automata,positioning,arrows}
\usepackage{young}
\usepackage[object=vectorian]{pgfornament} 
\usepackage{lipsum,tikz}
\usepackage{multirow}
\usepackage[OT2,T1]{fontenc}
\usepackage{mathtools}
\usepackage{ytableau}

\hypersetup{
	colorlinks=true, 
	linktoc=all, 
	linkcolor=blue} 

\numberwithin{equation}{section}

\theoremstyle{theorem}
\newtheorem{theorem}{Theorem}[section]
\newtheorem*{theorem*}{Theorem}

\newtheorem{corollary}[theorem]{Corollary}

\newtheorem{innercustomgeneric}{\customgenericname}
\providecommand{\customgenericname}{}
\newcommand{\newcustomtheorem}[2]{%
	\newenvironment{#1}[1]
	{%
		\renewcommand\customgenericname{#2}%
		\renewcommand\theinnercustomgeneric{##1}%
		\innercustomgeneric
	}
	{\endinnercustomgeneric}
}
\newcustomtheorem{ctheorem}{Theorem}
\newcustomtheorem{clemma}{Lemma}

\theoremstyle{definition}

\newtheorem*{example*}{Example}
\newtheorem*{examples*}{Examples}

\newtheorem*{remark*}{Remark}
\newtheorem*{remarks*}{Remarks}
\newtheorem*{note*}{Note}

\newtheoremstyle{named}{}{}{\itshape}{}{\bfseries}{.}{.5em}{\thmnote{#3} #1}
\theoremstyle{named}

\newtheoremstyle{customized}{}{}{\itshape}{}{\bfseries}{.}{.5em}{\thmnote{#3}}
\theoremstyle{customized}

\newcommand{\arxiv}[1]{\href{https://arxiv.org/abs/#1}{arXiv:#1}}
\newcommand{\doi}[1]{\href{https://dx.doi.org/#1}{DOI: #1}}

\DeclareMathAlphabet{\mydutchcal}{U}{dutchcal}{m}{n}
\newcommand{\dcC}{\mydutchcal{C}}

\newcommand{\qbinom}[2]{{\genfrac{[}{]}{0pt}{}{#1}{#2}}}
\newcommand{\tqbinom}[2]{{\textstyle\genfrac{[}{]}{0pt}{}{#1}{#2}}}

\newcommand{\GK}{\operatorname{GK}}

\newcommand{\bc}{\boldsymbol{c}}

\title[Cylindric partitions of rank three]{The bivariate generating function for cylindric partitions of rank three}

\author[S. Chern]{Shane Chern}
\address[S. Chern]{Fakult\"at f\"ur Mathematik, Universit\"at Wien, Oskar-Morgenstern-Platz 1, Wien 1090, Austria}
\email{chenxiaohang92@gmail.com, xiaohangc92@univie.ac.at}

\author[T. Nguyen]{Trieu Nguyen}
\address[T. Nguyen]{MathVision AI}
\email{trieu@mathvision.ai}

\author[C. Tran]{Chanh Tran}
\address[C. Tran]{MathVision AI}
\email{chanh@mathvision.ai}

\author[T. Wakhare]{Tanay Wakhare}
\address[T. Wakhare]{{}\textsuperscript{(1)}Department of Electrical Engineering and Computer Science, MIT, Cambridge, MA 02139, USA \newline \indent{}\textsuperscript{(2)}MathVision AI}
\email{tanay@mathvision.ai, twakhare@mit.edu}

\thanks{Mathematical contributors: S. Chern, T. Wakhare; Engineering contributors: T. Nguyen, C. Tran.}

\date{}



\begin{document}
	
\sloppy

\begin{abstract}
	Using a rational analogue of an $\mathrm{A}_2^{(1)}$ Rogers--Ramanujan type identity discovered by Warnaar, we establish sum-like expressions for a bivariate generating function for cylindric partitions of rank three. These expressions provide the finishing touch towards conjectures of Warnaar and Kanade--Russell.
\end{abstract}

\maketitle

\section{Introduction}

An \emph{integer partition} is a weakly decreasing sequence of positive integers known as \emph{parts}, and we allow the sequence to be empty, giving the \emph{empty partition}. The \emph{size} of a partition is the sum of all its parts, and the \emph{length} is the number of parts. We further tacitly assume that a part is zero if its index exceeds the length of this partition. Given an integer partition, its \emph{Young diagram} is obtained by depicting left-aligned boxes in rows such that the number of boxes in the $i$-th row equals the $i$-th part of the partition. We may generalize integer partitions in a two-dimensional way as \emph{plane partitions} by filling the Young diagrams with positive integers such that all rows and columns are weakly decreasing.

Inspired by the work of Burge~\cite{Bur1993}, Gessel and Krattenthaler~\cite{GK1997} introduced cylindric partitions as an affine analogue of plane partitions. These combinatorial objects have close connections to representation theory, especially to the principal characters of standard $\mathrm{A}_{r-1}^{(1)}$ modules; see, for example, \cite{FW2016,Tin2008}. Let $\bc := (c_0,c_1,\ldots,c_{r-1})$ be a sequence of nonnegative integers. We say a sequence $\Lambda := (\lambda^{(0)},\lambda^{(1)},\ldots,\lambda^{(r-1)})$ of integer partitions is a \emph{cylindric partition} of \emph{profile} $\bc$ and \emph{rank} $r$ if the inequalities
\begin{align*}
	\lambda_j^{(i)}\ge \begin{cases}
		\lambda_{j+c_{i+1}}^{(i+1)}, & i<r-1,\\[6pt]
		\lambda_{j+c_{0}}^{(0)}, & i=r-1,
	\end{cases}
\end{align*}
hold for all $j\ge 1$. We define the \emph{size} of $\Lambda$, written as $|\Lambda|$, as the overall sizes of the partitions $\lambda^{(0)},\ldots,\lambda^{(r-1)}$, and the \emph{maximum part} of $\Lambda$, written as $\max(\Lambda)$, as the maximum of all parts among these $r$ partitions. Let $\dcC_{\bc}$ be the set of cylindric partitions of profile $\bc$. We focus on the bivariate generating function
\begin{align*}
	\GK_{\bc}(z,q) := \sum_{\Lambda\in \dcC_{\bc}} z^{\max(\Lambda)} q^{|\Lambda|}.
\end{align*}
In particular, by the inherent cyclic symmetry of cylindric partitions, we have
\begin{align*}
	\GK_{(c_0,c_1,\cdots,c_{r-1})}(z,q) = \GK_{(c_1,\cdots,c_{r-1},c_0)}(z,q).
\end{align*}

It was known from the work of Borodin~\cite[Proposition~5.1]{Bor2007} that the univariate specialization of $\GK_{\bc}(z,q)$ at $z=1$ admits a product formula:
\begin{align*}
	&\GK_{(c_0,c_1,\ldots,c_{r-1})}(1,q)\\
	&\qquad = \frac{(q^{r+c_0+\cdots+c_{r-1}};q^{r+c_0+\cdots+c_{r-1}})_\infty^{r-1}}{(q;q)_\infty^r} \prod_{1\le i<j \le r} \theta(q^{j-i+c_i+\cdots+c_{j-1}};q^{r+c_0+\cdots+c_{r-1}}),
\end{align*}
where
\begin{align*}
	\theta(z;q) := (z;q)_\infty (q/z;q)_\infty,
\end{align*}
and
\begin{align*}
	(z)_\infty = (z;q)_\infty := \prod_{k\ge 0} (1-zq^k), \qquad\qquad (z)_n = (z;q)_n := \frac{(z;q)_\infty}{(zq^n;q)_\infty}.
\end{align*}
Notably, the infinite product in Borodin's formula first appeared in the work of Andrews, Schilling, and Warnaar~\cite[eq.~(5.2)]{ASW1999} in the context of the $\mathrm{W}_n$ algebra. Borodin's product formula has given the field of multiple Rogers--Ramanujan type identities a new impetus. Quite recently, these Rogers--Ramanujan type identities and the related cylindric partitions also emerged in an important geometric setting, arising from the Quot schemes associated with torus knot singularities~\cite{Huang2026,HJO2026}. In particular, for the curve $X^r = Y^s$ with $1<r<s$ coprime, the corresponding counting is closely tied to the enumeration of certain rank-$r$ cylindric partitions~\cite{HLO2026,HLPO2026,LO2026}.

For the moment, we restrict our attention to cylindric partitions of rank three. In \cite[Conjectures~2.5 and 2.11]{War2023}, Warnaar first conjectured bivariate Rogers--Ramanujan type summations for $\GK_{(k\pm 1,k,k)}(z,q)$ and $\GK_{(3k-s+1\pm 1,s-1,0)}(z,q)$. Kanade and Russell~\cite[Conjecture~5.1]{KR2023} expanded this list of conjectural expressions to all $\GK_{(c_0,c_1,c_2)}(z,q)$ with $c_1,c_2\le \lfloor\frac{c_0+c_1+c_2+1}{3}\rfloor$. Moreover, Kanade and Russell showed in \cite[Theorem~5.4]{KR2023} that using their conjectured cases as seeds and applying a recursive relation discovered by Corteel and Welsh~\cite[Proposition~3.1]{CW2019}, the bivariate generating function $\GK_{(c_0,c_1,c_2)}(z,q)$ for all possible profiles $(c_0,c_1,c_2)$ can be written as a linear combination of $z$-deformed Rogers--Ramanujan type multisums.

In \cite[Theorem~1.2]{War2025}, Warnaar gave a proof of a version of the Kanade--Russell seed conjecture. More precisely, he proved a Rogers--Ramanujan type identity~\cite[eq.~(1.5)]{War2025} whose sum side comes from the specialization $\GK_{(c_0,c_1,c_2)}(1,q)$ of Kanade and Russell's conjecture. Although Warnaar has offered all necessary ingredients, especially a rational analogue of an $\mathrm{A}_2^{(1)}$ Rogers--Ramanujan type identity~\cite[eq.~(5.14)]{War2025}, to derive the bivariate generating function $\GK_{(c_0,c_1,c_2)}(z,q)$, this task was not explicitly addressed in \cite{War2025}. Therefore, the purpose of the present work is to continue with Warnaar's method to formulate the desired $\GK_{(c_0,c_1,c_2)}(z,q)$ in the original conjecture of Kanade and Russell. 

\begin{theorem}[Kanade--Russell, Conjecture~5.1 in \cite{KR2023}]\label{th:KR-2.5}
	Let $(c_0,c_1,c_2)\in \mathbb{Z}_{\ge 0}^3$ be such that $c_1,c_2\le k$ where
	\begin{align*}
		c_0+c_1+c_2=3k+\epsilon, \qquad\qquad \epsilon\in \{-1,0,1\}.
	\end{align*}
	Define
	\begin{align*}
		\Omega^{(\epsilon)}(q) := \begin{cases}
			\qbinom{m_k+n_k}{m_k}_{q^3}\qbinom{m_k+n_k}{m_k}_q^{-1}, & \text{if $\epsilon = 0$},\\
			1, & \text{if $\epsilon = \pm 1$},
		\end{cases}
	\end{align*}
	where
	\begin{align*}
		\qbinom{n}{m}_q:=\begin{cases}
			\dfrac{(q;q)_n}{(q;q)_m(q;q)_{n-m}}, & \text{if $0\le m\le n$},\\[10pt]
			0, & \text{otherwise}.
		\end{cases}
	\end{align*}
	If $c_2=0$,
	\begin{align}\label{eq:GK-c2=0}
		&\GK_{(c_0,c_1,c_2)}(z,q)\notag\\
		&\qquad= \frac{(q)_\infty}{(zq)_\infty} \sum_{\substack{m_1,\ldots,m_k\ge 0\\n_1,\ldots,n_k\ge 0}} \frac{q^{\epsilon(\epsilon-1)m_kn_k+\sum\limits_{i=1}^k (m_i^2-m_in_i+n_i^2+m_i) +\sum\limits_{i=c_1+1}^k n_i}}{(q)_{m_1}(q)_{n_1}(q)_{m_k+n_k+1}}\notag\\
		&\qquad\quad\times \Omega^{(\epsilon)}(q)\qbinom{m_1}{m_2}_q\cdots\qbinom{m_{k-1}}{m_k}_q\qbinom{n_1}{n_2}_q\cdots\qbinom{n_{k-1}}{n_k}_q z^{n_1};
	\end{align}
	if $c_1=0$ and $c_2>0$,
	\begin{align}\label{eq:GK-c1=0&c2>0}
		&\GK_{(c_0,c_1,c_2)}(z,q)\notag\\
		&\qquad= \frac{(q)_\infty}{(zq)_\infty} \sum_{\substack{m_1,\ldots,m_k\ge 0\\n_1,\ldots,n_k\ge 0}} \frac{q^{\epsilon(\epsilon-1)m_kn_k+\sum\limits_{i=1}^k (m_i^2-m_in_i+n_i^2+n_i) +\sum\limits_{i=c_2+1}^k m_i}}{(q)_{m_1}(q)_{n_1}(q)_{m_k+n_k+1}}\notag\\
		&\qquad\quad\times \Omega^{(\epsilon)}(q)\qbinom{m_1}{m_2}_q\cdots\qbinom{m_{k-1}}{m_k}_q\qbinom{n_1}{n_2}_q\cdots\qbinom{n_{k-1}}{n_k}_q z^{n_1} \big(1-(1-z)q^{m_{c_2}+n_1+1}\big);
	\end{align}
	and if $c_1,c_2>0$,
	\begin{align}\label{eq:GK-c1>0&c2>0}
		&\GK_{(c_0,c_1,c_2)}(z,q)\notag\\
		&\qquad= \frac{(q)_\infty}{(zq)_\infty} \sum_{\substack{m_1,\ldots,m_k\ge 0\\n_1,\ldots,n_k\ge 0}} \frac{q^{\epsilon(\epsilon-1)m_kn_k+\sum\limits_{i=1}^k (m_i^2-m_in_i+n_i^2)+\sum\limits_{i=c_2+1}^k m_i +\sum\limits_{i=c_1+1}^k n_i}}{(q)_{m_1}(q)_{n_1}(q)_{m_k+n_k+1}}\notag\\
		&\qquad\quad\times (1-q^{m_{c_2}+n_{c_1}+1})\,\Omega^{(\epsilon)}(q)\qbinom{m_1}{m_2}_q\cdots\qbinom{m_{k-1}}{m_k}_q\qbinom{n_1}{n_2}_q\cdots\qbinom{n_{k-1}}{n_k}_q z^{n_1}.
	\end{align}
\end{theorem}

As corollaries, we also establish the alternative sum-like expressions conjectured by Warnaar~\cite{War2023} for $\GK_{(k\pm 1,k,k)}(z,q)$ and $\GK_{(3k-s+1\pm 1,s-1,0)}(z,q)$.

\begin{corollary}[Warnaar, Conjectures~2.5(a) and 2.11(a) in \cite{War2023}]\label{coro:Warnaar-a}
	We have
	\begin{align}\label{eq:Warnaar-a-(-)}
		&\GK_{(k-1,k,k)}(z,q)\notag\\
		&= \frac{1}{(zq)_\infty} \sum_{\substack{m_1,\ldots,m_{k-1}\ge 0\\n_1,\ldots,n_k\ge 0}} \frac{z^{n_1} q^{n_k^2+\sum\limits_{i=1}^{k-1} (m_i^2-m_in_i+n_i^2)}}{(q)_{n_1}}\notag\\
		&\quad\times \qbinom{n_1}{n_2}_q\cdots\qbinom{n_{k-1}}{n_k}_q\qbinom{n_1-n_2+m_2}{m_1}_q\cdots\qbinom{n_{k-2}-n_{k-1}+m_{k-1}}{m_{k-2}}_q\qbinom{n_{k-1}+n_{k}}{m_{k-1}}_q,
	\end{align}
	and
	\begin{align}\label{eq:Warnaar-a-(+)}
		&\GK_{(k+1,k,k)}(z,q)\notag\\
		&= \frac{1}{(zq)_\infty} \sum_{\substack{m_1,\ldots,m_k\ge 0\\n_1,\ldots,n_k\ge 0}} \frac{z^{n_1} q^{\sum\limits_{i=1}^{k} (m_i^2-m_in_i+n_i^2)}}{(q)_{n_1}}\notag\\
		&\quad\times \qbinom{n_1}{n_2}_q\cdots\qbinom{n_{k-1}}{n_k}_q\qbinom{n_1-n_2+m_2}{m_1}_q\cdots\qbinom{n_{k-1}-n_{k}+m_{k}}{m_{k-1}}_q\qbinom{2n_{k}}{m_{k}}_q.
	\end{align}
\end{corollary}

\begin{corollary}[Warnaar, Conjectures~2.5(b) and 2.11(b) in \cite{War2023}]\label{coro:Warnaar-b}
	We have, for $1\le s\le k+1$,
	\begin{align}\label{eq:Warnaar-b-(-)}
		&\GK_{(3k-s,s-1,0)}(z,q)\notag\\
		&= \frac{1}{(zq)_\infty} \sum_{\substack{m_1,\ldots,m_{k-1}\ge 0\\n_1,\ldots,n_k\ge 0}} \frac{z^{n_1} q^{n_k^2+\sum\limits_{i=1}^{k-1} (m_i^2-m_in_i+n_i^2+m_i)+\sum\limits_{i=s}^{k}n_i}}{(q)_{n_1}}\notag\\
		&\quad\times \qbinom{n_1}{n_2}_q\cdots\qbinom{n_{k-1}}{n_k}_q\qbinom{n_1-n_2+m_2}{m_1}_q\cdots\qbinom{n_{k-2}-n_{k-1}+m_{k-1}}{m_{k-2}}_q\qbinom{n_{k-1}+n_{k}}{m_{k-1}}_q,
	\end{align}
	and
	\begin{align}\label{eq:Warnaar-b-(+)}
		&\GK_{(3k+2-s,s-1,0)}(z,q)\notag\\
		&= \frac{1}{(zq)_\infty} \sum_{\substack{m_1,\ldots,m_k\ge 0\\n_1,\ldots,n_k\ge 0}} \frac{z^{n_1} q^{\sum\limits_{i=1}^{k} (m_i^2-m_in_i+n_i^2+m_i)+\sum\limits_{i=s}^{k}n_i}}{(q)_{n_1}}\notag\\
		&\quad\times \qbinom{n_1}{n_2}_q\cdots\qbinom{n_{k-1}}{n_k}_q\qbinom{n_1-n_2+m_2}{m_1}_q\cdots\qbinom{n_{k-1}-n_{k}+m_{k}}{m_{k-1}}_q\qbinom{2n_{k}}{m_{k}}_q.
	\end{align}
\end{corollary}

\section{Proof of Theorem~\ref{th:KR-2.5}}

For two integer partitions $\mu$ and $\nu$, we say $\nu\subseteq \mu$ if $\nu_j\le \mu_j$ holds for every index $j\ge 1$. The \emph{skew Young diagram} $\mu/\nu$ with $\nu\subseteq \mu$ is obtained by removing the Young diagram of $\nu$ from the upper-left corner of the diagram of $\mu$. A \emph{skew plane partition} on $\mathbb{Z}_{\ge 0}$ of shape $\mu/\nu$ can be viewed as a filling of the diagram $\mu/\nu$ with \emph{nonnegative} integers such that the rows and columns are both weakly decreasing. Since the rows are weakly decreasing, we may encode the filling as a multipartition $(\lambda^{(0)},\lambda^{(1)},\ldots)$ on $\mathbb{Z}_{\ge 0}$, in which the length of $\lambda^{(i)}$ is $\mu_{i+1}-\nu_{i+1}$ and we allow the parts to be zero even if their indices do not exceed the length.

Let $r$ be a fixed positive integer. Fix also partitions $\mu,\nu$ with $\nu\subseteq \mu$ and an integer $d$ with
\begin{align*}
	d\ge \max\{\mu_1-\mu_r,\nu_1-\nu_r\}.
\end{align*}
A multipartition $(\lambda^{(0)},\lambda^{(1)},\ldots,\lambda^{(r-1)})$ on $\mathbb{Z}_{\ge 0}$ is a cylindric partition on $\mathbb{Z}_{\ge 0}$ of \emph{shape} $\mu/\nu/d$ if it is the encoded multipartition of a certain skew plane partition of shape $\mu/\nu$ and simultaneously $(\lambda^{(r-1)},\lambda^{(0)},\lambda^{(1)},\ldots,\lambda^{(r-1)})$ is the encoded multipartition of a certain skew plane partition of shape $(d+\mu_r,\mu_1,\ldots,\mu_r)/(d+\nu_r,\nu_1,\ldots,\nu_r)$.

Now fix the profile $\bc := (c_0,c_1,\ldots,c_{r-1})$. Define
\begin{align*}
	\GK_{\bc,N}(q) := \sum_{\substack{\Lambda\in \dcC_{\bc}\\\max(\Lambda)\le N}} q^{|\Lambda|}.
\end{align*}
It is clear that
\begin{align*}
	\GK_{\bc}(z,q) = 1 + \sum_{N\ge 1} \big(\GK_{\bc,N}(q) - \GK_{\bc,N-1}(q)\big) z^N,
\end{align*}
so that
\begin{align}\label{eq:GK-z}
	\GK_{\bc}(z,q) = (1-z) \sum_{N\ge 0} \GK_{\bc,N}(q) z^N.
\end{align}

Let
\begin{align*}
	d := c_0 + c_1 + \cdots + c_{r-1},
\end{align*}
and
\begin{align*}
	\nu := (c_1 + \cdots + c_{r-1}, c_2 + \cdots + c_{r-1}, \ldots, c_{r-1}, 0).
\end{align*}
Let $(L^r)$ be the rectangular partition of length $r$ with all parts equal to $L$. We know from \cite[eq.~(3.11)]{War2023} that
\begin{align*}
	\GK_{\bc,N}(q) = \lim_{L\to \infty} \sum_{\substack{\text{cylindric partition $\Lambda$ on $\mathbb{Z}_{\ge 0}$}\\\text{of shape $(L^r)/\nu/d$ and rank $r$}\\\text{with $\max(\Lambda)\le N$}}} q^{|\Lambda|}.
\end{align*}
In the meantime, a classic result of Gessel and Krattenthaler~\cite[Theorem~2]{GK1997} (see also \cite[eq.~(3.9)]{War2023}) tells us that the generating function on the right-hand side of the above equals
\begin{align*}
	\sum_{\substack{k_1,\ldots,k_r\in \mathbb{Z}\\k_1+\cdots+k_r=0}} \underset{1\le i,j\le r}{\det} \left(q^{r(d+r)\binom{k_i}{2}+(d+r)ik_i+(\nu_j-j)(rk_i+i-j)} \tqbinom{N+L-\nu_j-dk_i}{L-\nu_j-(d+r)k_i+j-i}_q\right).
\end{align*}
Therefore,
\begin{align*}
	\GK_{\bc,N}(q) = \sum_{\substack{k_1,\ldots,k_r\in \mathbb{Z}\\k_1+\cdots+k_r=0}} \underset{1\le i,j\le r}{\det} \left(\frac{q^{r(d+r)\binom{k_i}{2}+(d+r)ik_i+(\nu_j-j)(rk_i+i-j)}}{(q)_{N+rk_i+i-j}}\right).
\end{align*}
Now making the substitutions $k_i \mapsto -l_i$ for every $i$ and writing
\begin{align*}
	\rho_j := \nu_j+r-j = c_j + \cdots + c_{r-1} + r - j
\end{align*}
for every $j$, we obtain
\begin{align*}
	&\GK_{\bc,N}(q)\\
	&\quad= \sum_{\substack{l_1,\ldots,l_r\in \mathbb{Z}\\l_1+\cdots+l_r=0}} q^{\sum\limits_{i=1}^r \big(r(d+r)\binom{-l_i}{2} - (d+r)il_i + r(rl_i-i)\big) - \sum\limits_{j=1}^r (\rho_j-r)j} \underset{1\le i,j\le r}{\det} \left(\frac{q^{-\rho_j(rl_i-i)}}{(q)_{N-rl_i+i-j}}\right),
\end{align*}
which, using the condition $l_1+\cdots+l_r=0$, further gives
\begin{align*}
	\GK_{\bc,N}(q) = \sum_{\substack{l_1,\ldots,l_r\in \mathbb{Z}\\l_1+\cdots+l_r=0}} q^{\sum\limits_{i=1}^r \big(r(d+r)\binom{l_i}{2} - (d+r)il_i - r\rho_il_i\big)} \underset{1\le i,j\le r}{\det} \left(\frac{q^{(\rho_i-\rho_j)(rl_i-i)}}{(q)_{N-rl_i+i-j}}\right).
\end{align*}
In particular,
\begin{align}\label{eq:GK_c-3}
	&\GK_{(c_0,c_1,c_2),N}(q)\notag\\
	&\quad = \sum_{\substack{l_1,l_2,l_3\in \mathbb{Z}\\l_1+l_2+l_3=0}} q^{(c_0+c_1+c_2+3)\sum\limits_{i=1}^3\big(3\binom{l_i}{2} - il_i\big)-3\big((c_1+c_2+2)l_1+(c_2+1)l_2\big)}\notag\\
	&\quad\quad\times\det\begin{pmatrix}
		\frac{1}{(q)_{N-3l_1}} & \frac{q^{(c_1+1)(3l_1-1)}}{(q)_{N-3l_1-1}} & \frac{q^{(c_1+c_2+2)(3l_1-1)}}{(q)_{N-3l_1-2}}\\[5pt]
		\frac{q^{-(c_1+1)(3l_2-2)}}{(q)_{N-3l_2+1}} & \frac{1}{(q)_{N-3l_2}} & \frac{q^{(c_2+1)(3l_2-2)}}{(q)_{N-3l_2-1}}\\[5pt]
		\frac{q^{-(c_1+c_2+2)(3l_3-3)}}{(q)_{N-3l_3+2}} & \frac{q^{-(c_2+1)(3l_3-3)}}{(q)_{N-3l_3+1}} & \frac{1}{(q)_{N-3l_3}}
	\end{pmatrix}.
\end{align}

Now let us define an auxiliary summation:
\begin{align*}
	S_{M,N} &:= \sum_{\substack{l_1,l_2,l_3\in \mathbb{Z}\\l_1+l_2+l_3=0}} q^{(c_0+c_1+c_2+3)\sum\limits_{i=1}^3\big(3\binom{l_i}{2} - il_i\big)-3\big((c_1+c_2+2)l_1+(c_2+1)l_2\big)}\\
	&\ \quad\times \det\begin{pmatrix}
		\frac{1}{(q)_{M+3l_1}(q)_{N-3l_1}} & \frac{q^{(c_1+1)(3l_1-1)}}{(q)_{M+3l_1+1}(q)_{N-3l_1-1}} & \frac{q^{(c_1+c_2+2)(3l_1-1)}}{(q)_{M+3l_1+2}(q)_{N-3l_1-2}}\\[5pt]
		\frac{q^{-(c_1+1)(3l_2-2)}}{(q)_{M+3l_2-1}(q)_{N-3l_2+1}} & \frac{1}{(q)_{M+3l_2}(q)_{N-3l_2}} & \frac{q^{(c_2+1)(3l_2-2)}}{(q)_{M+3l_2+1}(q)_{N-3l_2-1}}\\[5pt]
		\frac{q^{-(c_1+c_2+2)(3l_3-3)}}{(q)_{M+3l_3-2}(q)_{N-3l_3+2}} & \frac{q^{-(c_2+1)(3l_3-3)}}{(q)_{M+3l_3-1}(q)_{N-3l_3+1}} & \frac{1}{(q)_{M+3l_3}(q)_{N-3l_3}}
	\end{pmatrix}.
\end{align*}
It is clear from \eqref{eq:GK_c-3} that
\begin{align}\label{eq:GK-S}
	\GK_{(c_0,c_1,c_2),N}(q) = (q)_\infty^3 \lim_{M\to\infty} S_{M,N}.
\end{align}
Recall that
\begin{align*}
	c_0+c_1+c_2=3k+\epsilon, \qquad\qquad \epsilon\in \{-1,0,1\},
\end{align*}
and
\begin{align*}
	0\le c_1,c_2\le k.
\end{align*}
By a rational analogue of an $\mathrm{A}_2^{(1)}$ Rogers--Ramanujan type identity discovered by Warnaar~\cite[eq.~(5.14)]{War2025}, in which we choose $(a,b)=(c_1,c_2)$ and $\tau=\epsilon$, it is straightforward to compute that
\begin{align}\label{eq:S-formula}
	S_{M,N} &= \frac{1}{(q)_{M+N}}\sum_{\substack{0\le m_k\le\cdots\le m_1\le M\\0\le n_k\le\cdots\le n_1\le N}} \frac{q^{\epsilon(\epsilon-1)m_kn_k+\sum\limits_{i=1}^k (m_i^2-m_in_i+n_i^2)+\sum\limits_{i=c_2+1}^k m_i +\sum\limits_{i=c_1+1}^k n_i}}{(q)_{M-m_1}(q)_{N-n_1}(q)_{m_1}(q)_{n_1}(q)_{m_k+n_k+1}}\notag\\
	&\quad\times (1-q^{m_{c_2}+n_{c_1}+1})\,\Omega^{(\epsilon)}(q)\qbinom{m_1}{m_2}_q\cdots\qbinom{m_{k-1}}{m_k}_q\qbinom{n_1}{n_2}_q\cdots\qbinom{n_{k-1}}{n_k}_q,
\end{align}
where we put $m_0 := M$ and $n_0:= N$.

If $c_2=0$, we derive from \eqref{eq:GK-S} and \eqref{eq:S-formula} that
\begin{align*}
	&\GK_{(c_0,c_1,c_2),N}(q)\notag\\
	&\quad = (q)_\infty \sum_{\substack{0\le m_k\le\cdots\le m_1\\0\le n_k\le\cdots\le n_1\le N}} \frac{q^{\epsilon(\epsilon-1)m_kn_k+\sum\limits_{i=1}^k (m_i^2-m_in_i+n_i^2+m_i) +\sum\limits_{i=c_1+1}^k n_i}}{(q)_{N-n_1}(q)_{m_1}(q)_{n_1}(q)_{m_k+n_k+1}}\notag\\
	&\quad\quad\times \Omega^{(\epsilon)}(q)\qbinom{m_1}{m_2}_q\cdots\qbinom{m_{k-1}}{m_k}_q\qbinom{n_1}{n_2}_q\cdots\qbinom{n_{k-1}}{n_k}_q.
\end{align*}
Thus, in view of \eqref{eq:GK-z},
\begin{align*}
	&\GK_{(c_0,c_1,c_2)}(z,q)\\
	&\quad= (1-z) (q)_\infty \sum_{\substack{m_1,\ldots,m_k\ge 0\\n_1,\ldots,n_k\ge 0}} \frac{q^{\epsilon(\epsilon-1)m_kn_k+\sum\limits_{i=1}^k (m_i^2-m_in_i+n_i^2+m_i) +\sum\limits_{i=c_1+1}^k n_i}}{(q)_{m_1}(q)_{n_1}(q)_{m_k+n_k+1}}\notag\\
	&\quad\quad\times \Omega^{(\epsilon)}(q)\qbinom{m_1}{m_2}_q\cdots\qbinom{m_{k-1}}{m_k}_q\qbinom{n_1}{n_2}_q\cdots\qbinom{n_{k-1}}{n_k}_q \sum_{N\ge n_1} \frac{z^N}{(q)_{N-n_1}}.
\end{align*}
By Euler's first sum~\cite[eq.~(II.1)]{GR2004},
\begin{align*}
	\sum_{N\ge n_1} \frac{z^N}{(q)_{N-n_1}} = z^{n_1} \sum_{N\ge 0} \frac{z^N}{(q)_{N}} = \frac{z^{n_1}}{(z)_\infty}.
\end{align*}
The desired expression in \eqref{eq:GK-c2=0} then follows.

If $c_1=0$ and $c_2>0$, we have
\begin{align*}
	&\GK_{(c_0,c_1,c_2),N}(q)\notag\\
	&\quad = (q)_\infty \sum_{\substack{0\le m_k\le\cdots\le m_1\\0\le n_k\le\cdots\le n_1\le N}} \frac{q^{\epsilon(\epsilon-1)m_kn_k+\sum\limits_{i=1}^k (m_i^2-m_in_i+n_i^2+n_i) +\sum\limits_{i=c_2+1}^k m_i}}{(q)_{N-n_1}(q)_{m_1}(q)_{n_1}(q)_{m_k+n_k+1}}\notag\\
	&\quad\quad\times (1-q^{N+m_{c_2}+1})\,\Omega^{(\epsilon)}(q)\qbinom{m_1}{m_2}_q\cdots\qbinom{m_{k-1}}{m_k}_q\qbinom{n_1}{n_2}_q\cdots\qbinom{n_{k-1}}{n_k}_q.
\end{align*}
It follows that
\begin{align*}
	&\GK_{(c_0,c_1,c_2)}(z,q)\\
	&\quad= (1-z) (q)_\infty \sum_{\substack{m_1,\ldots,m_k\ge 0\\n_1,\ldots,n_k\ge 0}} \frac{q^{\epsilon(\epsilon-1)m_kn_k+\sum\limits_{i=1}^k (m_i^2-m_in_i+n_i^2+n_i) +\sum\limits_{i=c_2+1}^k m_i}}{(q)_{m_1}(q)_{n_1}(q)_{m_k+n_k+1}}\notag\\
	&\quad\quad\times \Omega^{(\epsilon)}(q)\qbinom{m_1}{m_2}_q\cdots\qbinom{m_{k-1}}{m_k}_q\qbinom{n_1}{n_2}_q\cdots\qbinom{n_{k-1}}{n_k}_q \sum_{N\ge n_1} \frac{z^N(1-q^{N+m_{c_2}+1})}{(q)_{N-n_1}},
\end{align*}
which, after using Euler's first sum twice, gives us \eqref{eq:GK-c1=0&c2>0}.

If $c_1,c_2>0$, we have
\begin{align*}
	&\GK_{(c_0,c_1,c_2),N}(q)\notag\\
	&\quad = (q)_\infty \sum_{\substack{0\le m_k\le\cdots\le m_1\\0\le n_k\le\cdots\le n_1\le N}} \frac{q^{\epsilon(\epsilon-1)m_kn_k+\sum\limits_{i=1}^k (m_i^2-m_in_i+n_i^2)+\sum\limits_{i=c_2+1}^k m_i +\sum\limits_{i=c_1+1}^k n_i}}{(q)_{N-n_1}(q)_{m_1}(q)_{n_1}(q)_{m_k+n_k+1}}\notag\\
	&\quad\quad\times (1-q^{m_{c_2}+n_{c_1}+1})\,\Omega^{(\epsilon)}(q)\qbinom{m_1}{m_2}_q\cdots\qbinom{m_{k-1}}{m_k}_q\qbinom{n_1}{n_2}_q\cdots\qbinom{n_{k-1}}{n_k}_q.
\end{align*}
Therefore,
\begin{align*}
	&\GK_{(c_0,c_1,c_2)}(z,q)\\
	&\quad= (1-z) (q)_\infty \sum_{\substack{m_1,\ldots,m_k\ge 0\\n_1,\ldots,n_k\ge 0}} \frac{q^{\epsilon(\epsilon-1)m_kn_k+\sum\limits_{i=1}^k (m_i^2-m_in_i+n_i^2)+\sum\limits_{i=c_2+1}^k m_i +\sum\limits_{i=c_1+1}^k n_i}}{(q)_{m_1}(q)_{n_1}(q)_{m_k+n_k+1}}\notag\\
	&\quad\quad\times (1-q^{m_{c_2}+n_{c_1}+1})\,\Omega^{(\epsilon)}(q)\qbinom{m_1}{m_2}_q\cdots\qbinom{m_{k-1}}{m_k}_q\qbinom{n_1}{n_2}_q\cdots\qbinom{n_{k-1}}{n_k}_q \sum_{N\ge n_1} \frac{z^N}{(q)_{N-n_1}},
\end{align*}
thereby producing \eqref{eq:GK-c1>0&c2>0}, again, in light of Euler's first sum.

\section{Proofs of Corollaries~\ref{coro:Warnaar-a} and \ref{coro:Warnaar-b}}

\subsection{Proof of Corollary~\ref{coro:Warnaar-a}}

We know from Theorem~\ref{th:KR-2.5} that for $\epsilon = \pm 1$,
\begin{align*}
	\GK_{(k+\epsilon,k,k)}(z,q) &= \frac{(q)_\infty}{(zq)_\infty} \sum_{\substack{m_1,\ldots,m_k\ge 0\\n_1,\ldots,n_k\ge 0}} \frac{z^{n_1} q^{\epsilon(\epsilon-1)m_kn_k+\sum\limits_{i=1}^k (m_i^2-m_in_i+n_i^2)}}{(q)_{m_1}(q)_{n_1}(q)_{m_k+n_k}}\notag\\
	&\quad\times\qbinom{m_1}{m_2}_q\cdots\qbinom{m_{k-1}}{m_k}_q\qbinom{n_1}{n_2}_q\cdots\qbinom{n_{k-1}}{n_k}_q.
\end{align*}
The relations \eqref{eq:Warnaar-a-(-)} and \eqref{eq:Warnaar-a-(+)} then follow from \cite[Lemma~7.1]{War2023}, in which we take $m_0,n_0\to\infty$.

\subsection{Proof of Corollary~\ref{coro:Warnaar-b}}

We again derive from Theorem~\ref{th:KR-2.5} that for $\epsilon = \pm 1$,
\begin{align*}
	&\GK_{(3k-s+1+\epsilon,s-1,0)}(z,q)\notag\\
	&\quad= \frac{(q)_\infty}{(zq)_\infty} \sum_{\substack{m_1,\ldots,m_k\ge 0\\n_1,\ldots,n_k\ge 0}} \frac{z^{n_1} q^{\epsilon(\epsilon-1)m_kn_k+\sum\limits_{i=1}^k (m_i^2-m_in_i+n_i^2+m_i)+\sum\limits_{i=s}^k n_i}}{(q)_{m_1}(q)_{n_1}(q)_{m_k+n_k+1}}\notag\\
	&\quad\quad\times\qbinom{m_1}{m_2}_q\cdots\qbinom{m_{k-1}}{m_k}_q\qbinom{n_1}{n_2}_q\cdots\qbinom{n_{k-1}}{n_k}_q.
\end{align*}
Now single out the summations over the $m$'s,
\begin{align*}
	&\GK_{(3k-s+1+\epsilon,s-1,0)}(z,q)\notag\\
	&\quad= \frac{(q)_\infty}{(zq)_\infty} \sum_{n_1,\ldots,n_k\ge 0} \frac{z^{n_1}q^{\sum\limits_{i=1}^k n_i^2+\sum\limits_{i=s}^k n_i}}{(q)_{n_1}}\qbinom{n_1}{n_2}_q\cdots\qbinom{n_{k-1}}{n_k}_q\\
	&\quad\quad\times \sum_{m_1,\ldots,m_k\ge 0} \frac{q^{\epsilon(\epsilon-1)m_kn_k+\sum\limits_{i=1}^k (m_i^2-m_in_i+m_i)}}{(q)_{m_1}(q)_{m_k+n_k+1}} \qbinom{m_1}{m_2}_q\cdots\qbinom{m_{k-1}}{m_k}_q.
\end{align*}
Using \cite[Lemma~7.2]{War2023}, we find that
\begin{align*}
	&\sum_{m_1,\ldots,m_k\ge 0} \frac{q^{2m_kn_k+\sum\limits_{i=1}^k (m_i^2-m_in_i+m_i)}}{(q)_{m_1}(q)_{m_k+n_k+1}}\qbinom{m_1}{m_2}_q\cdots\qbinom{m_{k-1}}{m_k}_q\\
	&\qquad\qquad\qquad\qquad = \frac{1}{(q)_\infty}\sum_{m_1,\ldots,m_{k-1}\ge 0} q^{\sum\limits_{i=1}^{k-1} (m_i^2-m_in_i+m_i)}\\
	&\qquad\qquad\qquad\qquad\quad\times\qbinom{n_1-n_2+m_2}{m_1}_q\cdots\qbinom{n_{k-2}-n_{k-1}+m_{k-1}}{m_{k-2}}_q\qbinom{n_{k-1}+n_{k}}{m_{k-1}}_q,
\end{align*}
and
\begin{align*}
	&\sum_{m_1,\ldots,m_k\ge 0} \frac{q^{\sum\limits_{i=1}^k (m_i^2-m_in_i+m_i)}}{(q)_{m_1}(q)_{m_k+n_k+1}} \qbinom{m_1}{m_2}_q\cdots\qbinom{m_{k-1}}{m_k}_q\\
	&\qquad\qquad\qquad\qquad = \frac{1}{(q)_\infty}\sum_{m_1,\ldots,m_k\ge 0} q^{\sum\limits_{i=1}^k (m_i^2-m_in_i+m_i)}\\
	&\qquad\qquad\qquad\qquad\quad\times \qbinom{n_1-n_2+m_2}{m_1}_q\cdots\qbinom{n_{k-1}-n_{k}+m_{k}}{m_{k-1}}_q\qbinom{2n_{k}}{m_{k}}_q,
\end{align*}
and therefore arrive at the desired relations \eqref{eq:Warnaar-b-(-)} and \eqref{eq:Warnaar-b-(+)}.

\subsection*{Statement on AI use}

T.W.~acknowledges the use of MathVision AI, an interactive research workspace, which pointed him towards Warnaar's identity~\cite[eq.~(5.14)]{War2025} and produced a proof of \eqref{eq:Warnaar-a-(-)}. MathVision was subsequently used by T.W.~to explore this connection. On the other hand, after becoming aware of Warnaar's identity, S.C.~proved the results in this work independently, without further consulting the AI outputs. AI tools were judiciously used throughout the research workflow of T.W., while the final draft was fully human-written by the two mathematical contributors (S.C. and T.W.).

\subsection*{Acknowledgements}

Shane Chern was supported by the FWF Austrian Science Fund (No.~10.55776/F1002).

\bibliographystyle{amsplain}

\end{document}